\documentclass[12pt]{amsart}
\usepackage{amssymb}
\usepackage{epsfig}
\usepackage{graphicx}
\numberwithin{equation}{section}

\input xy
\xyoption{all}

\advance\headheight 2pt
\calclayout
\allowdisplaybreaks[3]

\theoremstyle{plain}
\newtheorem{prop}{Proposition}

\newtheorem{theo}[prop]{Theorem}
\newtheorem{coro}[prop]{Corollary}

\newtheorem{lemm}[prop]{Lemma}

\theoremstyle{definition}
\newtheorem{defi}[prop]{Definition}

\newtheorem{rema}[prop]{Remark}

\newcommand{\bA}{\mathbb A}

\newcommand{\bG}{\mathbb G}
\newcommand{\bN}{\mathbb N}
\newcommand{\bP}{\mathbb P}
\newcommand{\bQ}{\mathbb Q}
\newcommand{\bZ}{\mathbb Z}

\newcommand{\cC}{\mathcal C}
\newcommand{\cD}{\mathcal D}
\newcommand{\cE}{\mathcal E}
\newcommand{\cF}{\mathcal F}
\newcommand{\cI}{\mathcal I}
\newcommand{\cL}{\mathcal L}
\newcommand{\cM}{\mathcal M}
\newcommand{\cN}{\mathcal N}
\newcommand{\cO}{\mathcal O}

\newcommand{\cT}{\mathcal T}

\newcommand{\cX}{\mathcal X}
\newcommand{\cY}{\mathcal Y}

\newcommand{\tcD}{\tilde{\mathcal D}}
\newcommand{\tcX}{\tilde{\mathcal X}}

\newcommand{\ra}{\rightarrow}
\newcommand{\da}{\downarrow}

\newcommand{\lra}{\longrightarrow}

\newcommand{\stab}{{\rm stab}}

\newcommand{\Def}{{\rm Def}}
\newcommand{\Hom}{{\rm Hom}}
\newcommand{\iHom}{{\mathcal H}\!{\mathit om}}

\newcommand{\Ext}{{\rm Ext}}
\newcommand{\Pic}{{\rm Pic}}

\makeatother
\makeatletter

\author{Brendan Hassett}
\address{Department of Mathematics \\
Rice University, MS 136 \\
Houston, TX 77251-1892}
\email{hassett@rice.edu}

\author{Yuri Tschinkel}
\address{
                Mathematisches Institut \\
                Bunsenstr. 3-5  \\
        37073 G\"ottingen \\
                Germany \\and\\
Courant Institute, NYU\\
251 Mercer Str. \\ New York,  NY 10012}
\email{tschinkel@cims.nyu.edu}

\title[Log Fano varieties over function fields]{Log Fano varieties over function fields
of curves}

\begin{document}
\date{\today}

\begin{abstract}
Consider a smooth log Fano variety over the function field of a curve.  
Suppose that the boundary has positive normal bundle.  Choose an integral
model over the curve.  Then integral
points are Zariski dense, after removing an explicit finite set of points
on the base curve.
\end{abstract}

\maketitle
\tableofcontents

\section{Introduction}
\label{sect:introduction}

Let $k$ be an algebraically closed field of characteristic zero,
$B$ a smooth projective curve over $k$ with function field $F=k(B)$.  

Our point of departure is the following theorem, combining work
of Graber-Harris-Starr and Koll\'ar-Miyaoka-Mori
\cite{GHS,KMM}:
Let $X$ be a smooth projective rationally-connected variety over $F$.  Then $X(F)$ is
Zariski dense in $X$.
One central example is
Fano varieties, i.e., varieties with 
ample anticanonical class, which
are known to be rationally connected (see \cite[V.2.13]{Kollar}).
In this context, it is not necessary to pass to a field extension
to get rational points.  

When $F$ is a number field, it may be necessary 
to pass to a finite extension to get rational points;
there exist Fano varieties over $\bQ$
without rational points.  
Moreover, even for Fano threefolds potential Zariski density, i.e., density 
after a finite extension of $F$, is unknown in general. For some positive
results in this direction see \cite{Harris-T}, \cite{Bogomolov-T} and \cite{Hassett}.

In this paper we study Zariski density of {\em integral} points. 
Consider pairs $(X,D)$ consisting of a variety $X$ and a divisor $D\subset X$, and 
fix integral models 
$\pi: (\cX,\cD) \ra B$  (see Section \ref{sect:prelim} for the definition). 
An $F$-rational point $s\in X\setminus D$ gives rise to a section
$s: B\ra \cX$ of $\pi$, meeting $\cD$ in finitely many points.
As we vary $s$, 
$$s^{-1}(\cD)=\pi(s(B)\cap \cD)\subset B$$
may vary as well.    
Fixing a finite set $S\subset B$, an $S$-integral point of $(\cX,\cD)$ is
an $F$-rational point of $X$ such that $s^{-1}(\cD)\subset S$ (as sets).  
\begin{theo}
Let $F$ be the function field of smooth projective curve $B/k$.  
Let $(X,D)$ be a pair consisting of a smooth projective variety $X$
and a smooth divisor $D\subset X$, defined over $F$.
Assume that $-(K_X+D)$ is ample and 
the first Chern class $c_1(\cN_{D/X})$ is effective 
and nonzero.  Choose
an integral model $(\cX,\cD) \ra B$.  Then there
is an explicit finite set $S\subset B$ such that
$S$-integral points of $(\cX,\cD)$ are Zariski dense.
\end{theo}
Theorem~\ref{theo:main} makes precise how $S$ is chosen.  

This is a partial converse to the function-field version of Vojta's conjectures:
Integral points are not Zariski dense when the log canonical class $K_X+D$
is ample (see \cite{vojta} for the number-field case, with connections
to value-distribution theory).  
Very few density results for integral points over number fields
are available, most of them in dimension two (see \cite{Silverman}, 
\cite{Beukers}, \cite{Hassett-T}).

\

\noindent 
{\bf Acknowledgments:}  We are grateful to Dan Abramovich for helpful
conversations on the deformation theory used in this article.
The first author appreciates the hospitality of the 
Mathematics Institute of the University of G\"ottingen.  
The first author was partially supported by National Science Foundation Grants
0554491 and 0134259 and an Alfred P. Sloan Research Fellowship.
The second author was partially supported by National Science Foundation 
Grants 0554280 and 0602333.

\section{Integral models and statements of results}
\label{sect:prelim}

\begin{defi}
A {\em pair}
$(X,D)$ consists of a smooth projective variety 
and a reduced effective divisor with normal crossings.  
\end{defi}

Let $B$ be a smooth projective curve defined over an algebraically closed field $k$ 
of characteristic zero and $F=k(B)$ its function field.  

\begin{defi}
Let $(X,D)$ be a pair defined over $F$.  
An {\em integral model} 
$$
\pi:(\cX,\cD) \ra B
$$ 
consists of a flat
proper morphism from a normal variety
$\pi_X:\cX\ra B$ with generic fiber $X$,
and a closed subscheme $\cD\subset \cX$ such that
$\pi_D:=\pi_X|_D:\cD\ra B$ is flat and has generic 
fiber $D$.
\end{defi}
We emphasize that $\cD$ has no irreducible components
contained in fibers of $\pi_X$.  

For many applications, the model is dictated by the 
specific circumstances.  Given an embedding of $(X,D)$ in projective space
there is a natural choice of model:  The properness
of the Hilbert scheme yields extensions of $X$ and $D$ to schemes
flat and projective over $B$.
Locally on $B$, these are obtained by
`clearing denominators' in the defining
equations of $X$ and $D$.  Normalizing if necessary, we obtain a model
of $(X,D)$.  

\begin{defi}
Let $S$ be a finite subset of $B$.  An {\em $S$-integral point of 
$(\cX,\cD)$} is a section $s:B\ra \cX$ such that
$s^{-1}(\cD)\subset S$ as sets.
\end{defi}
Thus if $D=\emptyset$ then integral points are just sections
of $\cX\ra B$, which are $F$-rational points of $X$.

The following proposition is straightforward:
\begin{prop} \label{prop:choice}
Let $(\cX_1,\cD_1)$ and $(\cX_2,\cD_2)$ be integral models of $(X,D)$.
Let $T\subset B$ denote the set over which the birational map
$$(\cX_1,\cD_1)\dashrightarrow (\cX_2,\cD_2)$$
fails to be an isomorphism.  
$S$-integral points
of $(\cX_1,\cD_1)$ are mapped to
$(S\cup T)$-integral points
of $(\cX_2,\cD_2)$.  
If $S$-integral points
of $(\cX_1,\cD_1)$ are Zariski dense then
$(S\cup T)$-integral points
of $(\cX_2,\cD_2)$ are Zariski dense.  
\end{prop}
In particular, if we allow ourselves to enlarge the 
set $S$ then
Zariski-density of integral points is independent of the model.

\begin{defi}
A point $b\in B$ is of {\em good reduction} 
if the fibers $\cX_b=\pi_X^{-1}(b)$ and $\cD_b=\pi_D^{-1}(b)$
are smooth.
\end{defi}

\begin{theo} \label{theo:main}
Let $(X,D)$ be a pair over $F=k(B)$ satisfying the following:
\begin{itemize}
\item{$D$ is smooth and rationally connected;}
\item{the normal bundle $\cN_{D/X}$ is effective
and nontrivial.}
\end{itemize}
Given a model $\pi:(\cX,\cD)\ra B$, let $S$
be a nonempty finite set of points in $B$ containing 
the image of the singular locus of $(\cX,\cD)$.  Then
$S$-integral points of $(\cX,\cD)$ are Zariski dense.
\end{theo}
Note however that we allow
points of bad reduction outside $S$.  
For example, let $\cX=\bP^2_{[x,y,z]} \times \bP^1_{[s,t]}$ and 
$$\cD=\{s(x^2+yz)+t(y^2+xz)=0 \}.$$
The model $(\cX,\cD)$ is smooth but $\cD_{[s,t]}$ is singular 
when $s^3+t^3=0$.

Let $K_X$ denote the canonical class of $X$ and
$K_X+D$ the log canonical class of $(X,D)$.  
The pair $(X,D)$ is {\em log Fano} if $-(K_X+D)$ is ample.
By adjunction
$$(K_X+D)|D=K_D$$
so $-K_D$ is ample.  Thus 
$D$ is Fano hence rationally connected \cite{KMM} \cite[V.2.13]{Kollar}.
\begin{coro}
Let $(X,D)$ be a log Fano variety over $F$ with $X$ and $D$ smooth.   
Assume that $\cN_{D/X}$ is effective and nontrivial.  For each integral
model and collection of places as specified in Theorem~\ref{theo:main},
integral points are Zariski dense.
\end{coro}

We discuss how Theorem~\ref{theo:main} can be reduced to the case
of nonsingular integral models:
\begin{defi}
A {\em good resolution} of an integral model is a birational
proper morphism from a pair
$$\rho: (\tcX,\tcD) \ra (\cX,\cD)$$
such that
\begin{itemize}
\item{
$\rho^{-1}(\cD)=\tcD$;}
\item{$\rho$ is an isomorphism over the open subset 
of $(\cX,\cD)$ where $\cX$ is smooth and $\cD$ is 
normal crossings.}
\end{itemize}
\end{defi}

\begin{rema}  

\

\begin{enumerate}
\item{
$\tcD$ may very well have components contained
in fibers over $B$, so $(\tcX,\tcD)$ is {\em not}
necessarily an integral model.}
\item{
The normality assumption guarantees that
for each $b\in B$ and each irreducible component
of $\cX_b$, the total space $\cX$ is smooth
at the generic point of that component.  
In particular, $\rho$ is an
isomorphism over a dense open subset of each fiber.}
\end{enumerate}
\end{rema}

Let $(\tcX,\tcD)\ra (\cX,\cD)$ be a good resolution,
$S\subset B$ a finite set containing the images of the 
singularities of $\cX$ and $\cD$,
and $\tcD^{\circ}$ the union of the components of $\tcD$
dominating $B$.  We have:
\begin{itemize}
\item{$\tcD^{\circ}$ is normal crossings;}
\item{images of $S$-integral points of $(\tcX,\tcD^{\circ})$
under $\rho$ are $S$-integral points of $(\cX,\cD)$.}
\end{itemize}
We have a bijection
$$\rho:\tcX\setminus \tcD \ra \cX \setminus \cD, $$
so $S$-integral points of $(\cX,\cD)$ correspond to sections
$$\{\tilde{s}:B\ra \tcX:\tilde{s}^{-1}(\tcD)\subset S \}.$$
Since the fibral components of $\tcD$ lie over $S$,
$S$-integral points of $(\cX,\cD)$ are equal to 
$S$-integral points of $(\tcX,\tcD^{\circ})$.

This analysis reduces Theorem~\ref{theo:main} to:
\begin{theo}[Smooth case] \label{theo:second}
Retain the assumptions of Theorem~\ref{theo:main} and assume
in addition that $\cX$ and $\cD$ are nonsingular.
Then for any nonempty $S\subset B$ the $S$-integral
points in $(\cX,\cD)$ are Zariski dense.  
\end{theo}

\section{Atiyah classes and free curves}
We work over an algebraically closed field $k$.  

Let $C$ be a smooth projective variety with tangent sheaf $\cT_C$.
Its deformation space is denoted $\Def(C)$ and first-order
deformations are given by $H^1(C,\cT_C)$.  
Let $L$ be an invertible sheaf on $C$ and
$p:L\ra C$ the line bundle defined by the same cocycle.  
The deformation space of $L$ is denoted $\Def(L)$ 
and first-order deformations are given by $H^1(C,\cO_C)$.  

Let $\Def(C,L)$ denote deformations of both $C$ and $L$.  
Taking $\bG_m$-invariants of the tangent-bundle exact sequence
$$0\ra \cT_{L/C} \ra \cT_L \ra p^*\cT_C \ra 0$$
we obtain the Atiyah extension \cite{Atiyah}
\begin{equation} \label{eqn:Atiyah1}
0 \ra \cO_C \ra \cE_{C,L} \ra \cT_C \ra 0.
\end{equation}
This is classified by an element $\lambda \in \Ext^1(\cT_C,\cO_C)=H^1(C,\Omega^1_C)$,
which (up to sign) agrees with the Chern class $c_1(L)$ \cite[pp. 196]{Atiyah}.
First-order deformations of $(C,L)$ are given by $H^1(C,\cE_{C,L})$ 
\cite[pp. 241]{IllusieI}.  The homomorphisms in the long exact sequence
$$H^1(C,\cO_C) \ra H^1(C,\cE_{C,L}) \ra H^1(C,\cT_C)$$
are the differentials of natural morphisms of deformation spaces
$$\Def(L) \ra \Def(C,L) \ra \Def(C).$$

Consider the case where $C$ is a reduced projective scheme,
perhaps with singularities.  First-order
deformations are given by $\Ext^1(\Omega^1_C,\cO_C)$.  Working directly with
K\"ahler differentials rather than tangent bundles, we obtain a dual version
of the Atiyah extension \cite[pp. 241]{IllusieI}
\begin{equation} \label{eqn:Atiyah2}
0 \ra \Omega^1_C \ra \cF_{C,L} \ra \cO_C \ra 0,
\end{equation}
again classified by an element $\lambda = \pm c_1(L) \in H^1(C,\Omega^1_C).$
The long exact sequence
$$\Ext^1(\cO_C,\cO_C) \ra \Ext^1(\cF_{C,L},\cO_C) \ra \Ext^1(\Omega^1_C,\cO_C)$$
gives the differentials of
$$\Def(L) \ra \Def(C,L) \ra \Def(C).$$

Now let $C$ be a nodal curve embedded in a smooth projective variety $Y$;
the component of the Hilbert scheme parametrizing deformations of $C$ in $Y$ 
is denoted $\Def(C\subset Y)$.  
First-order deformations of $C\subset Y$ correspond to
$$\Hom(\cI_C/\cI_C^2,\cO_C)=\Gamma(C,\cN_{C/Y});$$ 
here $\cI_C$ is the ideal sheaf 
and $\cN_{C/Y}=\iHom(\cI_C/\cI_C^2,\cO_C)$ the normal sheaf. 
Let 
$\Def(C \ra Y)$
denote deformations of the {\em morphism}
$C\ra Y$;  
first-order deformations correspond to
$\Hom(\Omega^1_Y|C,\cO_C)=\Gamma(C,T_Y|C)$ (see \cite[I.2]{Kollar}).

Fix a line bundle $L$ on $Y$.  Consider the morphisms
$$\begin{array}{rcl}
\mu:\Def(C\subset Y) & \ra & \Def(C,L|C) \\
  \{ C' \subset Y \} & \mapsto & (C',L|C')
\end{array}
$$
and
$$\begin{array}{rcl}
\nu:\Def(C\ra Y) & \ra & \Def(L|C) \\
  \{f:C'\ra  Y \} & \mapsto & f^*L.
\end{array}
$$
The same deformation space parametrizes fibers of both
$\mu$ and $\nu$:
\begin{equation} \label{eqn:defspace}
\begin{array}{rcl}
\Def((C,L|C)\ra (Y,L))&=&\{(f,M,\alpha):f \in \Def(C\ra Y), \\
	& & M \in \Def(L|C), 
    				 \alpha:M\stackrel{\sim}{\ra} f^*L \}.
\end{array}
\end{equation}

This has tangent space $\Hom(\cF_{Y,L}|C,\cO_C)=\Gamma(C,\cE_{Y,L}|C)$ and 
obstruction space $\Ext^1(\cF_{Y,L}|C,\cO_C)=H^1(C,\cE_{Y,L}|C)$.
Indeed, we have the diagram
\begin{equation} \label{eqn:Fdiagram}
\begin{array}{rcccccccl}
  &     &   0        &     &     0        &      &       &     &   \\
  &     &  \da       &     &    \da       &      &       &     &   \\
  &     & \cI_C/\cI_C^2 &=& \cI_C/\cI_C^2 &      &       &     &   \\
  &     &  \da       &     &    \da       &      &       &     &   \\
0 & \ra & \Omega^1_Y|C & \ra & \cF_{Y,L}|C & \ra & \cO_C & \ra & 0 \\
  &     &  \da       &     &    \da       &      &   ||  &     &   \\
0 & \ra & \Omega_C^1 & \ra & \cF_{C,L|C} & \ra & \cO_C & \ra & 0 \\
  &     &  \da       &     &    \da       &      &       &     &   \\
  &     &   0        &     &     0        &      &       &     &   
\end{array}
\end{equation}
where the middle row is obtained by restricting the dual Atiyah extension of
$(Y,L)$ to $C$.  Applying $\Hom(-,\cO_C)$ to this diagram, we obtain 
differentials between the various deformation spaces we have introduced.
The long exact sequence arising from the second row is
$$0 \lra \Gamma(\cO_C)\lra \Gamma(\cE_{Y,L}|C) \lra \Gamma(T_Y|C) 
\stackrel{d\nu}{\lra} H^1(\cO_C);$$
the second column yields
$$0 \lra \Gamma(\cE_{C,L|C}) \lra
\Gamma(\cE_{Y,L}|C)  \lra
\Gamma(\cN_{C/Y}) \stackrel{d\mu}{\lra}
H^1(\cE_{C,L|C}).
$$

A morphism is guaranteed to be smooth when its fibers are
unobstructed (cf. \cite[I.2.17.2]{Kollar}), thus we have
\begin{prop}
\label{prop:musmooth}
Let $C$ be a nodal curve embedded in a smooth projective variety
$Y$.  If  $H^1(C,\cE_{Y,L}|C)=0$ then 
$\mu$
is smooth at $C\subset Y$ and $\nu$ is smooth at $C\ra Y$.  
\end{prop}

\begin{defi}Let $Y$ be a smooth projective variety with
line bundle $L$ and $C$ a nodal projective curve.  
A nonconstant morphism $f:C\ra Y$ is 
{\em $L$-free} if for each $q\in C$
$$H^1(C,f^*\cE_{Y,L} \otimes \cI_q)=0.$$
It is {\em $L$-very free} if for each subscheme
$\Sigma \subset C$ of length two 
$$H^1(C,f^*\cE_{Y,L} \otimes \cI_{\Sigma})=0.$$
\end{defi}
Any $L$-free (resp. very free) morphism is free (resp. very free)
as $\cT_Y$ is a quotient of $\cE_{Y,L}$.

We now assume that $k$ is of characteristic zero.

\begin{prop} \label{prop:Lfree}
Let $Y$ be a smooth rationally connected projective variety,
$L$ a line bundle on $Y$, and $y \in Y$.  
Then $Y$ admits an $L$-free morphism
$f:\bP^1 \ra Y$ with image containing $y$.  
If $L$ is effective and nontrivial
then $f$ can be chosen to be $L$-very free.  
\end{prop}
\begin{proof}
There exists a very free morphism $g:\bP^1\ra Y$
\cite[IV.3.9.4]{Kollar};  moreover, given any
finite collection of points $y_1,\ldots,y_m \in Y$, 
we may assume the image of $g$ contains these points.

We have the extension
$$0 \ra \cO_{\bP^1} \ra g^*\cE_{Y,L} \ra g^*\cT_Y \ra 0$$
where $g^*\cT_Y$ is ample.  
It follows that each summand of $g^*\cE_{Y,L}$ is
nonnegative, which implies $L$-freeness.

Now assume $H$ is an effective
nonzero divisor corresponding to $L$.  
Choose $g$ such that its image contains $y$, some point
$y'$ in the support of $H$, and some point $y''$
not in the support of $H$.  In particular, the image
is not contained in any component of $H$.  It
follows that $g^*L$ has positive degree.  

If $\cO_{\bP^1}$ were a summand of $g^*\cE_{Y,L}$ then we 
would have 
$$g^*\cE_{Y,L} \simeq \cO_{\bP^1} \oplus g^*\cT_Y,$$
i.e., the Atiyah extension would split after pull-back.  
The extension induced by Diagram~\ref{eqn:Fdiagram}
$$0 \ra \cO_{\bP^1} \ra \cE_{\bP^1,g^*L} \ra \cT_{\bP^1} \ra 0$$
would split as well.  However, this extension is classified by
$$\pm c_1(g^*L) \in H^1(\bP^1,\Omega^1_{\bP^1})=\Ext^1(\cT_{\bP^1},\cO_{\bP^1}),$$
which is nontrivial.  
\end{proof}

A {\em comb with broken teeth} is a nodal projective curve 
$$C=B \cup T_1 \cup  \ldots \cup T_r$$
where $B$ is smooth and each $T_i$ is a tree of smooth 
rational curves meeting $B$
in a point $q_i \in B$.  Let $\sigma:C\ra B$ denote
the morphism which is the identity on $B$ and 
which contracts each $T_i$ to $q_i$.
There is a natural
stabilization morphism
$$\stab:\Def(C) \ra \Def(B)$$
constructed as follows:  Pick points 
$$b_1,\ldots,b_n \in B \setminus \{q_1,\ldots,q_r\}$$
with $2\mathsf{g}(B)-2+n>0$, so that $(B,b_1,\ldots,b_n)$
is a stable pointed curve.  
The Knudsen-Mumford stabilization
of $(C,\sigma^{-1}(b_1),\ldots,\sigma^{-1}(b_n))$ is $(B,b_1,\ldots,b_n)$. 
Each deformation of $C$
arises from a deformation of $(C,\sigma^{-1}(b_1),\ldots,\sigma^{-1}(b_n))$.
Thus the stabilization morphism
$$\Def(C,\sigma^{-1}(b_1),\ldots,\sigma^{-1}(b_n)) \ra \overline{M}_{\mathsf{g}(B),n}$$
induces $\stab$.  

For each
$e\in H^2(C,\bZ)$ (resp. $d\in H^2(B,\bZ)$), let $\Pic^e(C)$
(resp. $\Pic^d(B)$) denote the corresponding component of the
Picard scheme.  We have a morphism
$$\begin{array}{rcl}
\sigma_*:\Pic(C) & \ra & \Pic(B) \\
    M   & \mapsto & M|B\otimes \cO_B(e_1q_1+\ldots+e_rq_r),
                        \quad e_i=\deg(M|T_i)
\end{array}
$$
mapping each component $\Pic^e(C)$ isomorphically onto the component
$\Pic^d(B)$ containing its image.  Thus we get a morphism
$$\stab':\Def(C,M) \ra \Def(B,\sigma_*M)$$
and a commutative diagram
$$\begin{array}{ccc}
\Def(C,M) & \stackrel{\stab'}{\ra} & \Def(B,\sigma_*M) \\
\da & & \da \\
\Def(C) & \stackrel{\stab}{\ra} & \Def(B).
\end{array}
$$ 

Consider the composition
\begin{equation}
\label{whatistau}
\tau:\Def(C\subset Y) \stackrel{\mu}{\ra} \Def(C,L|C) \stackrel{\stab'}{\ra} 
		\Def(B,\sigma_*(L|C)).
\end{equation}
A fiber of $\tau$ corresponds to deformations of $C$ which do not
affect the line bundle induced by push-forward to the stabilization.  
\begin{prop} \label{prop:tauconstant}
Let $B$ be a smooth projective curve embedded in a smooth
variety $Y$ and $L$ a line bundle on $Y$.  
Assume that $H^1(B,\cE_{Y,L}|B)=0$ and consider a comb
$$C=B\cup T_1 \cup \ldots \cup T_r$$
such that the $T_i$ are $L$-free curves on $Y$.  Then the fiber 
$\tau^{-1}(\tau(C))$ contains a smoothing of $C$.  
\end{prop}
\begin{proof}
Write $q_i=B\cap T_i$ so that 
$$H^1(T_i,\cE_{Y,L} \otimes \cO_{T_i}(-q_i))=0.$$
Our vanishing assumption and an induction on the number of
components imply (see \cite[II.7.5]{Kollar})
$H^1(C,\cE_{Y,L}|C)=0$. 

We describe a flat morphism $\varpi:\cC \ra \bA^r$ deforming $C$
to $B$ (see \cite[pp. 156]{Kollar}).  Consider the smooth 
codimension-two subvariety
$$Z=\cup_{i=1}^r (\{q_i\} \times \{t_i=0\}) \subset B\times \bA^r$$
and the blow-up 
$$\sigma:\cC:=\mathrm{Bl}_Z(B\times \bA^r) \ra B\times \bA^r$$ 
with exceptional
divisors $E_1,\ldots,E_r$.  The composed morphism 
$$\varpi:\cC \ra B\times \bA^r \ra \bA^r$$
is still flat with $\varpi^{-1}(0)=C$;  {\em every}
fiber of $\varpi$ is a comb with handle $B$ and the
blow-down map is the stabilization contraction relative to $\bA^r$.  
We introduce a line bundle
on this family:  Consider
$$L'=L|B\otimes \cO_B(e_1q_1+\ldots+e_rq_r)$$ 
where $e_i=L\cdot T_i$ and write
$$\cM=(\pi_B\circ \sigma)^*L' \otimes \cO_{\cC}(-e_1E_1-\ldots-e_rE_r).$$
This is chosen so that 
$$\cM|\varpi^{-1}(0)=\cM|C=L|C$$
and $\sigma_*\cM=\pi_B^*L'$.  

We state a relative version of the 
deformation space~\ref{eqn:defspace}:
Consider morphisms
$$f:\cC \ra Y\times \bA^r$$
over $\bA^r$ admitting 
an isomorphism $\alpha:\cM \ra f^*(\pi_Y^*L)$.  
This is represented by a scheme
$$\Def((\cC,\cM) \ra(Y\times \bA^r, \pi_Y^*L)) \ra \bA^r$$
over $\bA^r$.  
The vanishing $H^1(C,\cE_{Y,L}|C)=0$ shows this problem
is unobstructed over $0\in \bA^r$ and thus the
deformation space is smooth over a neighborhood of $0 \in \bA^r$.  
In particular, it contains a smoothing of $C$ to $B$.  

By construction, the image of 
$$\Def((\cC,\cM)\ra (Y\times \bA^r,\pi_Y^*L))\ra \Def(C\subset Y)$$
is contained in the fiber of $\tau$.  
\end{proof}

The vanishing condition of Proposition~\ref{prop:tauconstant}
also guarantees that $\mu$ is smooth (see Proposition~\ref{prop:musmooth}).  
We indicate how to achieve this in practice
\begin{prop} \label{prop:getfree}
Let $Y$ be a smooth projective rationally connected variety,
$L$ an effective nontrivial line bundle, and $B$ a smooth
proper curve embedded in $Y$.  Then there exists a comb
$$C=B\cup T_1 \cup \ldots \cup T_r$$
such that $C$ deforms to a smooth $L$-free curve.  
In particular, 
$\tau$ is dominant at $(C\subset Y)\in \Def(C\subset Y)$.  
\end{prop}
\begin{proof}  Proposition~\ref{prop:Lfree}
gives $L$-very free rational curves through each point of $Y$.  
We use these to construct a comb with handle $B$ and $n\gg 0$
$L$-very free teeth $T_1,\ldots,T_n$.  
The Hard Smoothing technique of \cite[II.7.10]{Kollar}
implies that a subcomb
$$B \cup T_{i_1} \cup \ldots \cup T_{i_r}$$
deforms to a smooth $L$-free curve. 

Recall that $\tau=\stab' \circ \mu$ (see (\ref{whatistau})) and
$\mu$ is smooth near $C$.  Since $\stab'$ is birational, 
$\tau$ is dominant. 
\end{proof}

\section{Relative Atiyah classes and free curves}
In this section we work over an algebraically closed field
of characteristic zero.
See \cite{IllusieI}
for general background on relative obstruction theory.  

Let $B$ be a smooth projective curve.
Fix a proper nodal curve
$C$ over $B$.  Let $\Def(C/B)$ denote deformations of $C$ over $B$;
first-order deformations 
are parametrized by $\Ext^1(\Omega^1_{C/B},\cO_C)$.
If $L$ is a line bundle on $C$ then 
we have the relative Atiyah extension
$$0 \ra \Omega^1_{C/B} \ra \cF_{C,L/B} \ra \cO_C \ra 0.$$
Let $\Def(C,L/B)$ denote deformations
of $C$ and $L$ over $B$;
$\Ext^1(\cF_{C,L/B},\cO_C)$ parametrizes the first-order deformations.

Let $\pi:\cY \ra B$ be a nonconstant
proper morphism from a smooth variety.  
Given an embedding $C\subset \cY$ over $B$, we consider
the deformation space
$$\Def(C\ra \cY/B)$$
parametrizing deformations of the map $C\ra \cY$ over $B$.  This has
tangent space $\Gamma(C,\cT_{\cY/B}|C)$ and
obstruction space $H^1(C,\cT_{\cY/B}|C)$.  
Deformations of $C$ as a subscheme of $\cY$ are the same
as deformations of $C$ as a subscheme of $\cY$ over $B$, i.e.,
$$\Def(C\subset \cY/B)=\Def(C\subset \cY).$$

Let $\cL$ be an invertible sheaf on $\cY$
such that the restriction to the generic fiber of $\pi$ is
effective and nonzero.  
We have a relative version of the Atiyah extension
$$0 \ra \cO_{\cY} \ra \cE_{\cY,\cL/B} \ra \cT_{\cY/B} \ra 0,$$
classified by the image of the ordinary Atiyah class
under the restriction
$$\Ext^1(\cT_{\cY},\cO_Y) \ra \Ext^1(\cT_{\cY/B},\cO_Y).$$
Our effectivity assumption guarantees the extension is
not split.  

Consider the morphisms
$$
\mu_{/B}:\Def(C\subset Y)  \ra  \Def(C,L|C/B) 
$$
and
$$
\nu_{/B}:\Def(C\ra Y/B)  \ra  \Def(L|C), 
$$
i.e., the relative versions of $\mu$ and $\nu$ defined above.  
Their fibers are given by the relative version of (\ref{eqn:defspace}):
\begin{equation}
\label{eqn:reldefspace}
\begin{array}{rcl}
\Def((C,\cL|C)\ra (\cY,\cL)/B)&\!\!\!=\!\!\!&
	\{(f,M,\alpha):f \in \Def(C\ra \cY/B), \\
        & & M \in \Def(\cL|C),
                                 \alpha:M\stackrel{\sim}{\ra} f^*\cL \},
\end{array}
\end{equation}
with tangent space $\Gamma(C,\cE_{\cY,\cL/B}|C)$ and
obstruction space $H^1(C,\cE_{\cY,\cL/B}|C)$.  Here
$\cE_{\cY,\cL/B}=\cF_{\cY,\cL/B}^*$ where $\cF_{\cY,\cL/B}$ is
defined by the relative analog 
of Diagram~\ref{eqn:Fdiagram}:
\begin{equation} \label{eqn:Fdiagramrel}
\begin{array}{rcccccccl}
  &     &   0        &     &     0        &      &       &     &   \\
  &     &  \da       &     &    \da       &      &       &     &   \\
  &     & \cI_C/\cI_C^2 &=& \cI_C/\cI_C^2 &      &       &     &   \\
  &     &  \da       &     &    \da       &      &       &     &   \\
0 & \ra & \Omega^1_{\cY/B}|C & \ra & \cF_{\cY,\cL/B}|C & \ra & \cO_C & \ra & 0 \\
  &     &  \da       &     &    \da       &      &   ||  &     &   \\
0 & \ra & \Omega_{C/B}^1 & \ra & \cF_{C,\cL|C/B} & \ra & \cO_C & \ra & 0 \\
  &     &  \da       &     &    \da       &      &       &     &   \\
  &     &   0        &     &     0        &      &       &     &   
\end{array}
\end{equation}
Just as before, we obtain:
\begin{prop}
\label{prop:musmoothrel}
Retain the notation introduced above and assume that 
$H^1(C,\cE_{Y,L/B}|C)=0$.
Then the morphisms $\mu_{/B}$ and $\nu_{/B}$
are smooth at $C\subset \cY$.  
\end{prop}

\begin{defi} 
A nonconstant morphism from a nodal curve $f:C\ra \cY$ is 
{\em free over $B$} if for each $q\in C$
$$H^1(C,f^*\cT_{\cY/B} \otimes \cI_q)=0.$$
It is {\em very free over $B$} if for each subscheme
$\Sigma \subset C$ of length two 
$$H^1(C,f^*\cT_{\cY/B} \otimes \cI_{\Sigma})=0.$$
It is {\em $\cL$-free} or {\em $\cL$-very free over $B$}
if the analogous conditions hold for $\cE_{\cY,\cL/B}$.
\end{defi}

From now on, we will assume that
the generic fiber of $\pi:\cY \ra B$ is rationally connected.
\begin{prop}  \label{prop:Lfreerel}  Choose $b\in B$ such that
$\cY_b=\pi^{-1}(b)$ is smooth, and $y\in \cY_b$.  Then there exists
a morphism $f:\bP^1 \ra \cY_b$ that is $\cL$-very free over $B$,
with $y\in f(\bP^1)$.  
\end{prop}
\begin{proof} Our blanket assumptions imply $\cL|\cY_b$ is nontrivial
in every fiber of $\pi$.  (Nonzero divisors cannot specialize to zero in
smooth fibers.)   Proposition~\ref{prop:Lfree} then gives an 
$\cL|\cY_b$-very free curve in $\cY_b$.  Since 
$\cT_{\cY/B}|\cY_b=\cT_{\cY_b}$, we conclude this
curve is $\cL$-very free over $B$.
\end{proof}

We will require
the following relative version of Proposition~\ref{prop:getfree}:
\begin{prop} \label{prop:getfreerel}
Let $B'\subset \cY$ be a smooth projective curve,
not contained in a singular fiber of $\pi$.  
Then there exists a comb
$$C=B'\cup T_1 \cup \ldots \cup T_r$$
with teeth contained in smooth fibers of $\pi$,
such that $C$ deforms to a smooth curve that is
$\cL$-free over $B$.  
\end{prop}
\begin{proof} (cf. \cite[pp.63]{GHS})
If $B'$ is contained in a smooth fiber $\cY_b$ of $\pi$, we can apply 
Proposition~\ref{prop:getfree} with $Y=\cY_b$ and $B=B'$.  
Generally, the same proof applies:  Choose $\cL$-very free 
curves
$T_1,\ldots,T_n$ in smooth fibers $\cY_{b_1},\ldots,\cY_{b_n}$,
meeting $B'$ at points $q_1,\ldots,q_n$.  Deformations of the comb
$$B'\cup_{q_1} T_1 \cup_{q_2} T_2 \ldots \cup_{q_n} T_n$$
may be obstructed, but some subcomb
$$B'\cup T_{i_1} \cup\ldots \cup T_{i_r}$$
will smooth.  Again, the method of \cite[II.7.10]{Kollar} guarantees
that the generic such smoothing is $\cL$-free.  
\end{proof}

Let $C=B'\cup T_1 \cup \ldots \cup T_r$ be a comb with handle
$B'$ and broken teeth $T_1,\ldots,T_r$ contained in fibers
$\cY_{b_1},\ldots,\cY_{b_r}$.  
In analogy to (\ref{whatistau}) we define
\begin{equation}
\tau_{/B}:\Def(C\subset \cY) \stackrel{\mu_{/B}}{\ra} 
	\Def(C,\cL|C/B) \stackrel{\stab'}{\ra} 
		\Def(B',\sigma_*(\cL|C)/B),
\end{equation}
where $\sigma:C\ra B'$ contracts the teeth of the comb
and $\stab'$ is the stabilization introduced previously.
We can now state
the relative formulation of Proposition~\ref{prop:tauconstant}:
\begin{prop} \label{prop:tauconstantrel}
Retain the notation introduced above and consider a comb
$$C=B'\cup T_1 \cup \ldots \cup T_r$$
such that each $T_i$ is a $\cL$-free curve over $B$ in a smooth
fiber $\cY_{b_i}$.  
If $H^1(B',\cE_{\cY,\cL/B}|B')=0$ then the fiber
$\tau_{/B}^{-1}(\tau_{/B}(C))$ contains a smoothing of $C$.  
\end{prop}
\begin{rema}
If the $T_i$ are $\cL$-very free then we may apply 
Proposition~\ref{prop:getfreerel} to show the smoothing of $C$
is $\cL$-free over $B$. 
\end{rema}

Our main application is to sections of rationally-connected fibrations:
\begin{theo} \label{theo:prepare}
Let $B$ be a smooth projective curve, $\pi:\cY \ra B$ a proper
morphism from a smooth variety with rationally connected
generic fiber, and $\cL$ an invertible sheaf on $\cY$
restricting to a nontrivial effective divisor on the generic
fiber of $\pi$.  Fix an integer $N\gg 0$.

For each invertible sheaf $M\in \Pic^N(B)$, there exists a section
$s:B\ra \cY$ such that
$s^*\cL=M$ and
$s$ is $\cL$-free over $B$.  
In particular, the sheaves 
$$s^*\cE_{\cY,\cL/B} \text{ and } s^*\cT_{\cY/B}=\cN_{s(B)/\cY} $$
are both globally generated with no higher
cohomology.  
\end{theo}
In particular, Proposition~\ref{prop:musmooth} shows that the morphism
$$\begin{array}{rcl}
\mu:\Def(s(B)\subset \cY) & \ra & \Pic^N(B) \\
	s_t(B)  & \mapsto & s_t^*\cL 
\end{array}
$$
is smooth at $s(B)$.  
\begin{proof}[Proof of \ref{theo:prepare}]
The Graber-Harris-Starr Theorem \cite{GHS}
gives a section $s_1:B \ra \cY$.
The exact sequence
$$0 \ra \cT_{\cY/B} \ra \cT_{\cY} \ra \pi^*\cT_B \ra 0$$
induces
$$0 \ra s_1^*\cT_{\cY/B} \ra s_1^*\cT_{\cY} \ra \cT_B \ra 0$$
which is split by the differential $ds_1:\cT_B \ra s_1^*\cT_{\cY}$.  
Thus we have
$$s_1^*\cT_{\cY} = s_1^*\cT_{\cY/B} \oplus \cT_B$$
and the first term coincides with the normal bundle
$\cN_{s_1(B)/\cY}$.  

Proposition~\ref{prop:getfreerel} yields a section $s_2:B\ra \cY$
that is $\cL$-free over $B$ so that
$$H^1(B,s_2^*\cE_{\cY,\cL/B})=0.$$

To complete the proof, we apply Proposition~\ref{prop:tauconstantrel}
to produce a smoothing $s:B\ra \cY$ of a comb 
$$C=s_2(B) \cup_{b_1} T_1 \ldots \cup_{b_r}T_r.$$
However, it is necessary to relate the points of attachment
to the precise value of $s^*\cL$.  
Choose $e$ sufficiently large
so that for each smooth fiber $\cY_b$ and every point $y\in \cY_b$,
there exists an $\cL$-very free curve $T\subset \cY_b$ passing through $y$
with $\cL\cdot T=e$.  
(In Proposition~\ref{prop:Lfree}, we explained how to ensure
that $T$ intersects $\cL$ positively.)  
We therefore may assume 
$$e=T_1\cdot \cL=T_2\cdot \cL=\ldots =T_r\cdot \cL>0$$
so that
$$\sigma_*(\cL|C)=(s^*_2\cL)(eb_1+\ldots+eb_r).$$
Recall that $\tau_{/B}$ was defined so that 
for deformations $C'\subset \cY$
in the fiber $\tau_{/B}$ containing $C\subset \cY$,
the divisor class $\sigma_*(\cL|C')\in \Pic(B)$ remains constant.
Thus we have
$$s^*\cL=s^*_2\cL(eb_1+\ldots+eb_r).$$

Let $U\subset B$ denote the locus where $\cY_b$ is smooth
and contains a $\cL$-very free curve $T$ of degree $e$.  
It remains to verify the following prime
avoidance result, which governs the precise value of $N$:
\begin{lemm}
Let $B$ be a smooth projective curve,
$U\subset B$ a dense open subset, and $e\in \bN$.  
Fix a line bundle $\Lambda$ on $B$ of degree
$\ell$, $r\ge 2\mathsf{g}(B)+1$, and $N=er+\ell$.  For
any $M \in \Pic^N(B)$ there
exist distinct points
$b_1,\ldots,b_r \in U$ so that
$$M\simeq \Lambda(e(b_1+\ldots+b_r)).$$
\end{lemm}
\begin{proof} This is an elementary application of Riemann-Roch.
Choose an $e$th root of $M\otimes \Lambda^{-1}$, 
i.e., a line bundle $A$ with $A^{\otimes e}\otimes \Lambda=M$.
Any line bundle of degree $r$ on $B$ is very ample so consider
the embedding 
$$\phi_A:B \hookrightarrow \bP^{r-\mathsf{g}(B)}.$$
The divisors with any support along $B\setminus U$
form a finite union of hyperplanes in the linear system $|A|$.
The divisors admitting points of multiplicity $>1$ form
a proper subvariety of $\Delta \subset |A|$ by the Bertini Theorem.
Any divisor in the complement of the hyperplanes and $\Delta$
can be expressed in the form $b_1+\ldots+b_r$ with
the $b_i$ distinct in $U$.
\end{proof}
This concludes the proof of Theorem~\ref{theo:prepare}.
\end{proof}

\section{Smooth Case}
In this section, we prove Theorem~\ref{theo:second};
take $S=\{p\}$ for some $p \in B$.  

Apply Theorem~\ref{theo:prepare} to $\cY=\cD$, 
$\cL=\cO_{\cD}(\cD)=\cN_{\cD/\cX}$,
and $M=\cO_B(Np)$ for some suitable $N\gg 0$.  We obtain a section
$s:B\ra \cD$ with the following properties:
\begin{itemize}
\item{$s^*\cD=\cN_{\cD/\cX}|s(B) \simeq \cO_B(Np)$;}
\item{$\cE_{\cD,\cO(\cD)/B}|s(B)$ is globally generated with no
higher cohomology.}
\end{itemize}
Thus Proposition~\ref{prop:musmoothrel} guarantees 
$$\mu_{/B}:\Def(s(B) \subset \cD) \ra \Def(s(B),\cO_{\cD}(\cD)|s(B)/B)
		\simeq \Pic^N(B)$$
is smooth.  Let
$$\Def((B,M) \ra (\cD,\cL) /B)=\Def((B,\cO_B(Np)) \ra (\cD,\cO_{\cD}(\cD))/B)$$
denote the fiber, i.e.,
deformations $s_t:B\ra \cD$ such that 
$s_t^*\cD=\cO_B(Np)$.  

Consider the corresponding deformation space for $\cX$
$$\Def((B,\cO_B(Np)) \ra (\cX,\cO_{\cX}(\cD))/B);$$
its obstruction theory is
governed by the sheaf $\cE_{\cX,\cO(\cD)/B}|s(B)$.  
Consider the extensions defining $\cE_{\cD,\cO(\cD)/B}$ and $\cE_{\cX,\cO(\cD)/B}$,
restricted to $s(B)$:
$$\begin{array}{rcccccccl}
  &     &       &     &          0                 &     &       0           &     &   \\
  &     &       &     &     \downarrow             &     & \downarrow        &     &   \\
0 & \ra & \cO_B & \ra & \cE_{\cD,\cO(\cD)/B}|s(B) & \ra & \cN_{s(B)/\cD}    & \ra & 0 \\
  &     &   ||  &     &     \downarrow             &     & \downarrow        &     &   \\
0 & \ra & \cO_B & \ra & \cE_{\cX,\cO(\cD)/B}|s(B) & \ra & \cN_{s(B)/\cX}    & \ra & 0 \\
  &     &       &     &     \downarrow             &     & \downarrow        &     &   \\
  &     &       &     &    \cN_{\cD/\cX}|s(B)      &  =  &\cN_{\cD/\cX}|s(B) &     &   \\ 
  &     &       &     &     \downarrow             &     & \downarrow        &     &   \\
  &     &       &     &          0                 &     &       0           &     &   
\end{array}
$$ 
Since the terms in the bottom row are isomorphic to $\cO_B(Np)$, which has no
higher cohomology, we deduce that $\cE_{\cX,\cO(\cD)/B}|s(B)$ has no higher
cohomology.  In particular, 
$$\Def((B,\cO_B(Np)) \ra (\cX,\cO_{\cX}(\cD))/B)$$
is unobstructed and smooth.

The inclusion of
$\cD$ in $\cX$ induces an embedding
\begin{equation} \label{eqn:complicated}
\Def((B,\cO_B(Np)) \ra (\cD,\cO_{\cD}(\cD))/B) \hookrightarrow
\Def((B,\cO_B(Np)) \ra (\cX,\cO_{\cX}(\cD))/B).
\end{equation}
The image is precisely the indeterminacy of the rational map
$$\begin{array}{rcl}
G:\Def((B,\cO_B(Np)) \ra (\cX,\cO_{\cX}(\cD))/B)
& \dashrightarrow & 
		\bP(\Gamma(B,\cO_B(Np))) \\
     s_t(B) & \mapsto & s_t^*\cD.
\end{array}
$$
$S$-integral points are sections $s_t$
mapping to elements in $\Gamma(B,\cO_B(\cD))$ vanishing at $p$ to
maximal order $N$.  Thus we are interested in elements of
$G^{-1}[\Gamma(B,\cO_B)]$, where $\Gamma(B,\cO_B)\subset \Gamma(B,\cO_B(Np)$
corresponds to the constant functions, i.e., the image of the map 
on global sections induced by the inclusion of sheaves
$$\cO_B \hookrightarrow \cO_B(Np).$$

The indeterminacy of $G$ is resolved by blowing up the
subscheme (\ref{eqn:complicated}).
The stalk of its normal bundle 
at $s(B)$ is canonically isomorphic to $\Gamma(\cO_{\cX}(\cD)|s(B))$.  
In particular, the proper transform of $G^{-1}[\Gamma(B,\cO_B)]$ 
meets the exceptional fiber
over $s(B)$ at the point 
$$[\Gamma(B,\cO_B)] \in \bP(\Gamma(\cO_{\cX}(\cD)|s(B)))\simeq
\bP(\Gamma(B,\cO_B(Np))).$$
Thus $s(B)$ deforms to $s_t(B)\in G^{-1}[\Gamma(B,\cO_B)]$, 
corresponding to an $S$-integral point.

The sections thus produced are Zariski dense in $\cX$.  Indeed, our construction
produces sections passing through the generic point of $\cD$ that deform
out of $\cD$ to the generic point of $\cX$.

\bibliographystyle{plain}
\bibliography{LFD}

\begin{thebibliography}{10}

\bibitem{Atiyah}
Michael~F. Atiyah.
\newblock Complex analytic connections in fibre bundles.
\newblock {\em Trans. Amer. Math. Soc.}, 85:181--207, 1957.

\bibitem{Beukers}
Frits Beukers.
\newblock Integral points on cubic surfaces.
\newblock In {\em Number theory (Ottawa, ON, 1996)}, volume~19 of {\em CRM
  Proc. Lecture Notes}, pages 25--33. Amer. Math. Soc., Providence, RI, 1999.

\bibitem{Bogomolov-T}
Fedor Bogomolov and Yuri Tschinkel.
\newblock On the density of rational points on elliptic fibrations.
\newblock {\em J. Reine Angew. Math.}, 511:87--93, 1999.

\bibitem{GHS}
Tom Graber, Joe Harris, and Jason Starr.
\newblock Families of rationally connected varieties.
\newblock {\em J. Amer. Math. Soc.}, 16(1):57--67 (electronic), 2003.

\bibitem{Harris-T}
Joe Harris and Yuri Tschinkel.
\newblock Rational points on quartics.
\newblock {\em Duke Math. J.}, 104(3):477--500, 2000.

\bibitem{Hassett}
Brendan Hassett.
\newblock Potential density of rational points on algebraic varieties.
\newblock In {\em Higher dimensional varieties and rational points (Budapest,
  2001)}, volume~12 of {\em Bolyai Soc. Math. Stud.}, pages 223--282. Springer,
  Berlin, 2003.

\bibitem{Hassett-T}
Brendan Hassett and Yuri Tschinkel.
\newblock Density of integral points on algebraic varieties.
\newblock In {\em Rational points on algebraic varieties}, volume 199 of {\em
  Progr. Math.}, pages 169--197. Birkh\"auser, Basel, 2001.

\bibitem{IllusieI}
Luc Illusie.
\newblock {\em Complexe cotangent et d\'eformations. {I}}.
\newblock Springer-Verlag, Berlin, 1971.
\newblock Lecture Notes in Mathematics, Vol. 239.

\bibitem{Kollar}
J{\'a}nos Koll{\'a}r.
\newblock {\em Rational curves on algebraic varieties}, volume~32 of {\em
  Ergebnisse der Mathematik und ihrer Grenzgebiete. 3. Folge. A Series of
  Modern Surveys in Mathematics}.
\newblock Springer-Verlag, Berlin, 1996.

\bibitem{KMM}
J{\'a}nos Koll{\'a}r, Yoichi Miyaoka, and Shigefumi Mori.
\newblock Rationally connected varieties.
\newblock {\em J. Algebraic Geom.}, 1(3):429--448, 1992.

\bibitem{Silverman}
Joseph~H. Silverman.
\newblock Integral points on curves and surfaces.
\newblock In {\em Number theory (Ulm, 1987)}, volume 1380 of {\em Lecture Notes
  in Math.}, pages 202--241. Springer, New York, 1989.

\bibitem{vojta}
Paul Vojta.
\newblock {\em Diophantine approximations and value distribution theory},
  volume 1239 of {\em Lecture Notes in Mathematics}.
\newblock Springer-Verlag, Berlin, 1987.

\end{thebibliography}

\end{document}